\documentclass{amsart}
\usepackage{amsthm,amsfonts,amsmath,amscd,amssymb,latexsym,epsfig}
\newcommand{\gl}{{\mathfrak g \mathfrak l}}

\renewcommand{\u}{{\mathfrak u}}

\newcommand{\g}{{\mathfrak g}}         

\newcommand{\cx}{{\mathbb C}}
\newcommand{\diag}{\operatorname{diag}}

\newcommand{\Ad}{\operatorname{Ad}}
\newcommand{\tr}{\operatorname{tr}}

\newcommand{\Lie}{\operatorname{Lie}}

\newcommand{\Rat}{\operatorname{Rat_{\underline{k}}\bigl(F(n)\bigr)}}

\newcommand{\LPP}{\operatorname{\sL\sP\sP}}

\newcommand{\rat}{\operatorname{Rat\bigl(F(n)\bigr)}}

\numberwithin{equation}{section}

\newtheorem{theorem}{Theorem}[section]

\newtheorem{theorema}{Theorem}

\newtheorem{lemma}[theorem]{Lemma}

\newtheorem{corollary}[theorem]{Corollary}

\newtheorem{proposition}[theorem]{Proposition}

\theoremstyle{remark}

\newtheorem{remark}[theorem]{Remark}

\newtheorem{definition}[theorem]{Definition}

\newtheorem{ack}{Acknowledgment}

\newcommand{\oH}{{\mathbb{H}}}

\newcommand{\oN}{{\mathbb{N}}}

\newcommand{\oP}{{\mathbb{P}}}

\newcommand{\oR}{{\mathbb{R}}}

\newcommand{\oZ}{{\mathbb{Z}}}

\newcommand{\sG}{{\mathcal{G}}}   

\newcommand{\sL}{{\mathcal{L}}}   
\newcommand{\sM}{{\mathcal{M}}}   

\newcommand{\sO}{{\mathcal{O}}}
\newcommand{\sP}{{\mathcal{P}}}

\newcommand{\sR}{{\mathcal{R}}}

\newcommand{\fE}{{\mathfrak{e}}}

\newcommand{\fG}{{\mathfrak{g}}}
\newcommand{\fH}{{\mathfrak{h}}}

\newcommand{\fL}{{\mathfrak{l}}}

\newcommand{\fS}{{\mathfrak{s}}}

\begin{document}

\title{Gelfand-Zeitlin actions and rational maps}
\author{Roger Bielawski \and Victor Pidstrygach}
\address{School of Mathematics\\ University of Leeds\\Leeds LS2 9JT, UK}
\address{Mathematisches Institut, Universit\"at G\"ottingen, G\"ottingen 37073, Germany} 

\thanks{The research of the first author is supported by the Alexander von Humboldt Foundation}


\begin{abstract}  We observe that the analogue of the Gelfand-Zeitlin action on $\gl(n,\cx)$, which exists on any symplectic manifold $M$ with an Hamiltonian action of $GL(n,\cx)$, has a natural interpretation as a residual action, after we identify $M$ with a symplectic quotient of $M\times \prod_{i=1}^n T^\ast GL(i,\cx)$. We also show that the Gelfand-Zeitlin actions on $T^\ast GL(n,\cx)$ and on the regular part of $\gl(n,\cx)$ can be identified with natural Hamiltonian actions on spaces of rational maps into full flag manifolds, while the Gelfand-Zeitlin action on the whole $\gl(n,\cx)$ corresponds to a natural action on a space of rational maps into the manifold of half-full flags in $\cx^{2n}$.
\end{abstract}

\maketitle

\thispagestyle{empty}

\section{Introduction}
The Lie algebra $\gl(n,\cx)$ of $n\times n$ complex matrices admits a rather remarkable action of the group $\cx^{n(n-1)/2}$, which arose from the Gelfand-Zeitlin theory \cite{GZ}, and has been studied in detail by Kostant and Wallach in the two papers \cite{KW1} and \cite{KW2}. Let us briefly recall the definition of this Gelfand-Zeitlin group. 
\par
For an $n\times n$ complex matrix $B$ and for every $m\leq n$ let $B_{(m)}$ denote the upper-left $m\times m$ minor of $B$.  For every $1\leq i\leq m\leq n-1$, define now an action of $\cx$ on $\gl(n,\cx)$ via
\begin{equation} z.B=\Ad\left(e^{z(B_{(m)})^{i-1}}\right)B.\label{e}\end{equation}
The remarkable fact (cf. \cite{KW1}) is that the actions for different $(m,i)$ commute, giving an action of an abelian group  $A\simeq \cx^{n(n-1)/2}$. The action \eqref{e} is obtained by integrating the Hamiltonian vector field,  with respect to the standard (complex) Poisson structure on $\gl(n,\cx)\simeq\gl(n,\cx)^\ast$, corresponding to the function 
\begin{equation}  p_{mi}(B)=\tr \bigl(B_{(m)}\bigr)^i.\label{GZ}\end{equation}
Thus, the statement about commutativity of different \eqref{e}, translates into Poisson-commutativity of different $p_{mi}$. One can add the Hamiltonians $p_{n1},\dots,p_{nn}$, which have trivial Hamiltonian vector fields, and consider  the map $p=(p_{11},\dots,p_{nn}):\gl(n,\cx)\rightarrow \cx^{n(n+1)/2}$. This map is known as the Gelfand-Zeitlin map \cite{GZ,KW1}, and we shall refer to the action of $A$ as the {\em Gelfand-Zeitlin action}.  A beautiful ``symplectic" fact \cite{GKL0} about this map is that, on a Zariski open subset of regular adjoint orbit, the Hamiltonians $p_{11},\dots,p_{n-1,n-1}$ and the orbits of $A$ provide the ``angle-action coordinates", i.e. Darboux coordinates.
\par
We have found the Gelfand-Zeitlin group and its action, as described above,  quite mysterious  and the  purpose of this article is to give it some other, more geometric interpretations. 
\par 
First of all, an analogue of the Gelfand-Zeitlin action exists on any (complex) symplectic or Poisson manifold $M$  with an Hamiltonian action of $GL(n,\cx)$. The action is defined in the same way, i.e. via \eqref{GZ}, where now $B$ is replaced by the moment $\mu(p)$ at every point $p$ of $M$. The simplest way to see that this gives an action of the abelian $\cx^{n(n+1)/2}$  (without using the commutativity of the action on  $\gl(n,\cx)^\ast$) is to observe that $M$ is isomorphic to the symplectic quotient of $M\times \prod_{i=1}^n T^\ast GL(i,\cx)$. The Gelfand-Zeitlin action on $M$ is then the residual action of a group $ \cx^{n(n+1)/2}=\prod_{i=1}^n \cx^i$ with each factor acting only on the corresponding factor of $\prod_{i=1}^n T^\ast GL(i,\cx)$.
\par
We then turn to study the Gelfand-Zeitlin action on the symplectic manifold 
 \begin{equation}V_n=\{(B,b)\in \gl(n,\cx)\times\cx^n; \enskip\text{$b$ is a cyclic vector for B}\}\simeq GL(n,\cx)\times \cx^n.\label{V_n}\end{equation}
On this manifold, unlike on $\gl(n,\cx)^\ast$, we have the effective action of the full $ \cx^{n(n+1)/2}$, i.e. the vector fields corresponding to  $p_{n1},\dots,p_{nn}$ are nontrivial and give rise to an action of a group $C_n\simeq\cx^n$.  Furthermore, the symplectic quotients of $V_n$ by $C_n$ are precisely the regular adjoint orbits. Many properties of the Gelfand-Zeitlin action on $V_n$, and, hence, on regular adjoint orbits, become  more transparent thanks to the following result:
\begin{theorema} There exists a canonical symplectomorphism, equivariant with respect to the Gelfand-Zeitlin action, between $V_n$ and the space $\sR_n$ of based rational maps of degree $(1,\dots,n)$ from $\oP^1$ to the manifold of full flags in $\cx^{n+1}$.\label{theorema}\end{theorema}  
To explain what ``equivariant" means here, we have to describe an Hamiltonian action of $ \cx^{n(n+1)/2}$ on $\sR_n$.  Any such a rational map has a canonically defined polar part, consisting of $n$ monic polynomials $q_i(z)$ with $\deg q_i=i$. The coefficients of these polynomials provide a family of commuting Hamiltonians, whose flows define an action of a group $A^\prime \simeq \cx^{n(n+1)/2}$. Our results in the first part of the paper show that the action of $A^\prime $ on $\sR_n$, and similar actions on other spaces of rational maps, are much more transparent than the Gelfand-Zeitlin action on $V_n$ or on $\gl(n,\cx)$. In particular, the action respects the natural stratification of spaces of rational maps and, hence, it can be locally decomposed into actions of lower-dimensional groups.
We then reinterpret and improve some of the results of Kostant and Wallach using the above theorem.
\par
We give a similar interpretation of the Poisson structure and the Gelfand-Zeitlin action on $T^\ast GL(n,\cx)$ in terms of rational maps to the manifold of full flags in $\cx^{2n}$. We also observe that the Poisson structure of the  full $\gl(n,\cx)$ (and not only of its regular part)  is isomorphic to a certain space of rational maps into the manifold of {\em half-full flags} in $\cx^{2n}$. This isomorphism is again equivariant with respect to the Gelfand-Zeitlin action.

\tableofcontents

\part{Rational maps}

\section{Rational maps into  full flag manifolds\label{rff}}

For a detailed discussion of the spaces of rational  maps into flag manifolds see  \cite{MM,BHMM} and, for a somewhat different point of view, \cite{Hu}.

Let $\{e_1,\dots,e_{n+1}\}$ be a linear basis of $\cx^{n+1}$.  Consider the following  standard flag of subspaces of $\cx^{n+1}$:
$$ E_{0}^+=\{0\},\enskip E_1^+=\langle e_1\rangle,\dots,\enskip E_{n+1}^+=\cx^{n+1},$$
which we call the standard flag, and the anti-standard flag:
$${E}_i^-=\langle e_{n+1},\dots,e_{n-i+2}\rangle.$$
Let $B^+\leq GL(n+1,\cx)$ denote the stabiliser of the standard flag, i.e. the subgroup of  upper-triangular matrices and denote by $F(n)$ the homogeneous manifold $GL(n+1,\cx)/B^+$, i.e. the manifold of  full flags.  Denote by $B^-$ the subgroup of lower-triangular matrices and by $N^+,N^-$ the respective unipotent subgroups. The group $N^-$ acts on $F(n)$ and its orbits are known as the {\em Bruhat cells}. The open cell coincides with $N^-$ itself and the remaining ones  are in 1-1 correspondence with elements of the symmetric group $S_{n+1}$. In particular, the codimension one cells $X_1,\dots,X_{n}$ correspond to simple negative root spaces and their closures freely generate $H^2\bigl(F(n),\oZ\bigr)$. 
\par
We are interested in rational (i.e. holomorphic) maps $f:\oP^1\rightarrow F(n)$. A topological invariant of such a map is its degree $\underline{k}=(k_1,\dots,k_n)$,   where $k_i$ is the intersection number of $f\bigl(\oP^1\bigr)$ with $\overline{X}_i$. We base our rational maps by requiring that their value at $\infty$ is the anti-standard flag, i.e. $f(\infty)=({E}_i^-)$.
We denote by $\rat$ the space of based rational maps and by $\Rat$ its subspace of maps of degree  $\underline{k}$.  It is well-known that $\Rat$ is a complex manifold of dimension $2(k_1+\dots+k_n)$.

\medskip

{\sl Notation.} For a multi-index $\underline{k}=(k_1,\dots,k_n)$, we write $|\underline{k}|=k_1+\dots+k_n$.

\subsection{Poles and principal parts} Let $f\in \Rat$ and let $z_1^i,\dots, z_{k_i}^i\in \oP^1$ be the points which map to $\overline{X}_i$.  Thus, an $f\in \Rat$ defines $n$ polynomials:
\begin{equation} q_i(z)=\prod_{j=1}^{k_i}\bigl(z-z_j^i\bigr),\enskip i=1,\dots,n.\label{q}\end{equation}
Observe that, because of the basing condition, none of the $z_j^i$ is $\infty$, and, hence, each $q_i$ has degree $k_i$. Thus, associating to a rational map its polar divisor (or, equivalently, the coefficients of the $q_i(z)$, $i=1,\dots,n$), gives the  {\em pole location map}
\begin{equation} \Pi:\Rat\rightarrow  S^{|\underline{k}|}(\cx)\simeq \cx^{|\underline{k}|}.\label{Pi} \end{equation}
A based rational map $f$ is determined by the location of its poles and its {\em  local principal part} at every pole.  First of all let us classify the poles of $f$:
\begin{definition} Let $f:\oP^1\rightarrow F(n)$ be a rational map and let $\underline{r}=(r_1,\dots,r_n)$ be a sequence of nonnegative integers. We say that a point $z\in\oP^1$ is a {\em pole of type $\underline{r}$} of $f$, if $f\bigl(\oP^1\bigr)$ meets each $\overline{X_i}$  with multiplicity $r_i$ at $f(z)$. \label{type}\end{definition}

We now recall, after \cite{BHMM}, the concept of local principal parts, due to Gravesen \cite{Gr}. Let $Z$ denote the divisor at infinity of $F(n)$, i.e. the complement of the open cell $N^-$. Let $\sM$ denote the corresponding sheaf of germs of meromorphic maps into $N^-$, i.e. the sheaf corresponding to the presheaf $\sM(U)=\sO\bigl(U,F(n)\bigr)-\sO\bigl(U,Z\bigr)$,  and let $\sO(N^-)$ be the sheaf of germs of holomorphic maps into $N^-$. The action of $N^-$ on $F(n)$ induces a free action of  $\sO(N^-)$ on $\sM$, which leaves fixed the location and the order of poles. We obtain an exact sequence of sheaves 
\begin{equation} 0\rightarrow \sO(N^-)\rightarrow \sM\rightarrow \sP\sP\rightarrow 0,\label{PP}\end{equation}
and the quotient sheaf $\sP\sP$ is called the {\em sheaf of principal parts}. A based rational map to $F(n)$ can be viewed as a global section of $\sP\sP$ \cite{Gr,BHMM}. It is then clear that such map $f$ is determined by its poles and nontrivial elements of the stalks $\sP\sP_{z_i}$ at each pole $z_i$. These stalks are called local principal parts  and depend only on the type of a pole and not on its location. Thus, for an $\underline{r}=(r_1,\dots,r_n)$, we denote by $\sL\sP\sP_{\underline {r}}$ the stalk of $\sP\sP$ at a pole of type $\underline{r}$. We can also describe $\sL\sP\sP_{\underline {r}}$ (somewhat tautologically) as follows:
\begin{proposition} $\sL\sP\sP_{\underline{r}}=\left\{g\in \text{\rm Rat}_{\underline{r}}\bigl(F(n)\bigr) ;\enskip\text{all poles of $g$ are equal to $0$}\right\}$.\label{LPP}\end{proposition}

It is easy to verify, as in Example 3.12 of \cite{BHMM} (see also Lemma 4.9 in that paper), that
\begin{equation} \sL\sP\sP_{(0,\dots,1,\dots,0)}=\cx^\ast,\label{C*}\end{equation}
and, more generally, 
\begin{equation}\text{\rm Rat}_{(0,\dots,k,\dots,0)}\bigl(F(n)\bigr)\simeq  \text{\rm Rat}_k\bigl(\oP^1\bigr). \label{Rat_k}\end{equation}
In other words, a rational map $f\in \Rat$, whose every pole has pure type, i.e. $f$ does not meet codimension two Bruhat cells, is given by a collection 
\begin{equation}\left(\frac{p_1(z)}{q_1(z)},\dots,\frac{p_n(z)}{q_n(z)}\right)\label{generic}\end{equation}
of rational functions of degrees $k_1,\dots,k_n$. We shall also need to know the local principal parts and spaces of rational maps for two nonzero degrees:
\begin{proposition} Let $i\neq j$ and let $\underline{k}=(k_1,\dots,k_n)$ be such that $k_i=1$, $k_j=1$, and $k_m=0$ for $m\neq i,j$. Then
$$\LPP_{\underline{k}}=\begin{cases} \cx^\ast\times \cx^\ast &\text{if $|i-j|\geq 2$}\\  \cx^\ast\times\{(u,w)\in \cx^2; uw=0\} &\text{if $|i-j|=1$}.\end{cases}$$
Similarly
$$\Rat=\begin{cases}\bigl(\cx \times \cx^\ast\bigr)^2 &\text{if $|i-j|\geq 2$}\\  \bigl(\cx\times \cx^\ast\bigr)\times \cx^2 &\text{if $|i-j|=1$},\end{cases}$$
where, in the case $j=i+1$, the poles $z_i,z_{i+1}$ are related to $(Y,p,u,w)\in\bigl(\cx\times \cx^\ast\bigr)\times \cx^2$ by:
\begin{equation} z_{i+1}+z_i=Y, \enskip z_{i+1}-z_i=uw.\label{sts}\end{equation}
\label{(1,1)}\end{proposition}
\begin{proof} This is essentially given in Example 3.12 in \cite{BHMM}. Alternatively, one can use Hurtubise's identification, recalled in the Appendix, of $\Rat$ with a moduli space of solutions to Nahm's equations. In particular, it follows easily from Theorem \ref{Hurt} (see also Proposition \ref{oooo}), that $\text{\rm Rat}_{(1,1)}$ is the symplectic quotient of $V=(\cx\times\cx^\ast)^2\times \cx^2$ with the obvious symplectic structure by $\cx^\ast$, where $t\in \cx^\ast$ acts by on $\bigl((z_1,X_1),(z_2,X_2),(u,w)\bigr)\in V$ by
$$ V\ni \bigl((z_1,X_1),(z_2,X_2),(u,w)\bigr)\mapsto \left ((z_1,t^{-1}X_1),(z_2,tX_2),(tu,t^{-1}w) \right).$$
The second equation in \eqref{sts} is then the moment map equation for this action, and the quotient can be identified with $\bigl(\cx\times \cx^\ast\bigr)\times \cx^2$.\end{proof}

 \subsection{Stratification\label{stratification}}
The key advantage of the ``poles-and-principal-parts" point of view is their additivity. That is, if $f,g\in \rat$  have no poles in common, then viewing $f$ and $g$ as configurations of principal parts, their union $f\cup g$ is also a configuration of principal parts, giving an element of $\rat$. Thus, as in \eqref{generic}, we can consider subsets of $\Rat$, where rational maps decompose into ``sums" (i.e. unions in the above description) of rational maps of fixed lower degrees. This leads to the following stratification of $\Rat$ (cf. section 6 in \cite{BHMM}):
\par
Let $\overline{m}=\bigl\{\underline{m}^1,\dots,\underline{m}^r\bigr\}$ be a collection of multi-indices satisfying
\begin{equation} \sum_{i=1}^r\underline{m}^i=\underline{k}.\end{equation}
\begin{definition}\cite{BHMM} We denote by $R_{\overline{m}}$  the subset of all elements of $\Rat$ that have $r$ poles at distinct points $z_1,\dots,z_r\in \cx$ such that the type of the pole at $z_i$ is $\underline{m}^i$.\label{stratum}\end{definition}
Denote by $S^{|\underline{m^1}|, ...,|\underline{m^r}|}(\cx)$ the image $\Pi\bigl(R_{\overline{m}}\bigr)$ of the pole location map \eqref{Pi} in the symmetric product $S^{|\underline{k}|}(\cx)$. The pole location map  restricts to a locally trivial fibration \begin{equation}\Pi_{\overline{m}}:R_{\overline{m}}\rightarrow S^{|\underline{m^1}|, ...,|\underline{m^r}|}(\cx),\label{Pi_m}\end{equation} with fibre
\begin{equation} \prod_{i=1}^r\LPP_{\underline{m}^i}.\label{prod}\end{equation}
We now adopt the following definitions.
\begin{definition} For the collections of multi-indices $\overline{m}=\bigl\{\underline{m}^1,\dots,\underline{m}^r\bigr\}$ one has the addition map of divisors 
\begin{equation} \prod _{i=1}^r S^{|\underline{m}^i|}(\cx)\rightarrow S^{|\underline{k}|}(\cx).\end{equation} 
We shall consider its restriction $t$ to the subset $V(\overline{m}) \subset \prod _{i=1}^r S^{|\underline{m}^i|}(\cx)$ consisting of non-intersecting divisors:\begin{equation} t: V(\overline{m}) \rightarrow S^{|\underline{k}|}(\cx).\label{add-divisors}\end{equation}
Let 
\begin{equation}\Pi:\prod_{i=1}^r \text{\rm Rat}_{\underline{m}^i}\bigl(F(n)\bigr)\rightarrow \prod _{i=1}^r S^{|\underline{m}^i|}(\cx)\label{coord}\end{equation}
denote the coordinate-wise pole location map and let
\begin{equation} \sR(\overline{m})=\Pi^{-1}\bigl(V(\overline{m})\bigr).\label{sR}\end{equation}
In other words, $(f_1,\dots,f_r)\in \sR(\overline{m})$ if no poles of $f_i$ coincide with any poles of $f_j$, when $i\neq j$.  \end{definition}
\par
The ``partial addition" map, defined at the beginning of subsection, gives  a map
\begin{equation} T:\sR(\overline{m})\rightarrow \Rat.\label{T}\end{equation}
We observe that the stratum $R_{\overline{m}}$ is contained in the image of $T$.
\par
The following fact is rather straightforward:
\begin{proposition} The following diagram commutes:
\begin{equation}\begin{CD} \sR(\overline{m}) @>{T}>> \Rat \\
  @V{\Pi}VV  @V{\Pi}VV  \\
V(\overline{m}) @>{t}>> S^{|\underline{k}|}(\cx).\end{CD}\label{CD}\end{equation}
Moreover, any section of $t$ over a subset $U$ of $S^{|\underline{k}|}(\cx)$ lifts uniquely to a section of $T$ over $\Pi^{-1}(U)$.
\hfill $\Box$\label{TT}\end{proposition}

\subsection{Symplectic structure\label{symplectic}}
We denote by $\Rat^{\rm p}$ the open dense subset of  $\Rat$, consisting of those $f$ having all  poles of pure type, i.e. those $f$, which can be given by \eqref{generic}, and by $\Rat^\circ $ the open stratum of $\Rat$, the subset where all poles have multiplicity $1$. We have  canonical embeddings
\begin{equation} \Rat^\circ \hookrightarrow  \Rat^{\rm p}\hookrightarrow\prod_{i=1}^n   \text{\rm Rat}_{k_i}\bigl(\oP^1\bigr).\label{embed}\end{equation}
The second embedding can be interpreted as saying that, for $\underline{m}^i=(0,\dots,k_i,\dots 0)$, the map $T$ of \eqref{T} is an embedding.
\par
Hurtubise \cite{Hu} has given a biholomorphism between  $\Rat$ and the moduli space of $SU(n+1)$-monopoles with maximal symmetry breaking and charges $(k_1,\dots,k_n)$. Since the latter space is naturally a hyperk\"ahler manifold, one obtains a canonical (holomorphic) symplectic structure on $\Rat$. This symplectic structure has been computed on $\Rat^\circ $ by the first author in his thesis \cite{thesis}. It is simply the pullback of the product symplectic structure on $\prod_{i=1}^n   \text{\rm Rat}_{k_i}\bigl(\oP^1\bigr)$ under the map \eqref{embed}, i.e.:
\begin{equation} \omega=\sum_{i=1}^n\sum_{j=1}^{k_i}\frac{dp_i\bigl(z_j^i\bigr)}{p_i\bigl(z_j^i\bigr)}\wedge dq_i\bigl(z_j^i\bigr),\label{sympl}\end{equation}
where $z_i^j$ are the poles. Let us show directly, without monopoles, that this $\omega$ extends to the whole $\Rat$.
\begin{proposition} The form $\omega$, given by \eqref{sympl}, extends to a globally defined holomorphic symplectic form on $\Rat$.\label{Sympl}
\end{proposition}
\begin{proof} The form $\omega$ is defined on the complement of the codimension one variety consisting of rational maps with multiple poles. We claim that it extends to the complement  $\Rat - V$ of the codimension two variety consisting of  maps with triple or higher order poles. Let $f\in \Rat - V$ and suppose that $z$ is a double pole of $f$. Then, we have two possibilities for the type, in the sense of Definition \ref{type}: $\underline{r}=(r_1,\dots,r_n)$ of $z$: (1) $r_i=2$ and $r_j=0$ for $j\neq i$, or (2) $r_i=1$, $r_j=1$, and $r_m=0$ for $m\neq i,j$. In the first case, a neighbourhood of $f$ can be identified with an open subset in $\prod_m\text{\rm Rat}_{k_m}\bigl(\oP^1\bigr)$ and, hence, $\omega$ is defined at $f$, since it is defined everywhere on $\prod_m\text{\rm Rat}_{k_m}\bigl(\oP^1\bigr)$. In case (2), a neighbourhood of $f$ can be identified with an open subset in
$$ \prod_{m\neq i,j}\text{\rm Rat}_{k_m}\bigl(\oP^1\bigr)\times \text{\rm Rat}_{k_{i}}\bigl(\oP^1\bigr)\times \text{\rm Rat}_{k_{j}}\bigl(\oP^1\bigr)\times \text{\rm Rat}_{1,1}\bigl(F(2)\bigr).$$
Proposition \ref{(1,1)} implies that the symplectic form \eqref{sympl} on $\text{\rm Rat}_{1,1}\bigl(F(2)\bigr)^\circ$ extends to $\text{\rm Rat}_{1,1}\bigl(F(2)\bigr)$ and, hence, we have found an extension to a neighbourhood of $f$. 
\par
Thus $\omega$ defines a holomorphic map on the complement of a codimension two subvariety in $\Lambda^2 T\Rat$ and therefore extends to all of $\Lambda^2 T\Rat$ by Hartog's theorem. The extension is bilinear and closed by continuity. It is also non-degenerate on a complement of a codimension two subvariety, hence non-degenerate everywhere.\end{proof}

We finish the section by observing that the definition of $\omega$ implies that the map $T$ of \eqref{T} intertwines the relevant symplectic structures:
\begin{proposition} Let $T$ be the map \eqref{T}. The pullback $T^\ast \omega$ of the symplectic form of $\Rat$ is the product symplectic form of  $\prod_{i=1}^r \text{\rm Rat}_{\underline{m}^i}\bigl(F(n)\bigr)$.\hfill $\Box$\label{TTT}\end{proposition}

\section{Hamiltonians and the group $A_{\protect\underline{k}}$\label{A(c)}}

The formula \eqref{sympl} implies that the elementary symmetric polynomials in the pole location coordinates $z_1^i,\dots,z_{k_i}^i$,  $i=1,\dots,n$, are all Poisson-commuting. These polynomials are simply coefficients of $q_i$ given by \eqref{q}, $i=1,\dots,n$. Let us write 
\begin{equation} q_i=z^{k_i}+\sum_{j=1}^{k_i}c_{ij}z^{k_i-j},\quad\text{and}\enskip c=(c_{11},\dots,c_{nk_n}).\label{c_{ij}}\end{equation}
Thus, $f\mapsto c$ is a mapping from $\Rat$ to $\cx^{|\underline{k}|}$. It is nothing other than the pole location map \eqref{Pi}, after identifying $S^{|\underline{k}|}(\cx)$ with $\cx^{|\underline{k}|}$.

We now consider the flows of Hamiltonian vector fields corresponding to various $c_{ij}$:
\begin{proposition} The Hamiltonian vector field  corresponding to a $c_{ij}$, $i=1,\dots, n$, $j=1,\dots,k_i$,  is globally integrable defining an action of $\cx$ on $\Rat$. This action stabilises any fiber of the map $c_{ij}:\Rat \rightarrow \cx$.\label{global}\end{proposition}
\begin{proof} As in the proof of Proposition \ref{Sympl}, we can show that this is true on the Zariski open subset $U$ of $\Rat$, where rational maps have at most double poles. Thus, for  a $c_{ij}$, we obtain a holomorphic map $\lambda:\cx \times U\rightarrow \Rat$, given by integrating the Hamiltonian vector field corresponding to $c_{ij}$. Since $U$ is the complement of a codimension two subvariety, the map $\lambda$ extends to all of $\cx \times\Rat$. The second part of the statement is obvious.
\end{proof}

Since various $c_{ij}$ Poisson-commute, the flows corresponding to all of them define  an abelian subgroup $A_{\underline{k}}$ of the group of symplectomorphisms of $\Rat$. Of course, $A_{\underline{k}}\simeq \cx^{|\underline{k}|}$. We shall now study the orbits of $A_{\underline{k}}$.

\subsection{Orbits and stratification} 
We observe that on $\Rat^\circ $ (and only there), the orbits of $A_{\underline{k}}$ are isomorphic to $\bigl(\cx^\ast\bigr)^{|{\underline{k}}|}$ and, for a point of $\Rat^\circ $, its $A_{\underline{k}}$-orbit coincides with the fibre of the pole location map \eqref{Pi}. In general, any fibre $\Pi^{-1}(c)$ is stabilised by the group $A_{\underline{k}}$, but it is no longer a single orbit. We shall make use of the following observation:
\begin{proposition} The map $T$ of \eqref{T} maps any orbit of $\prod_{i=1}^r A_{\underline{m^i}}$ biholomorphically to an orbit of $A_{\underline{k}}$. \label{TTTT}\end{proposition}
\begin{proof} This follows from the commutativity of diagram \eqref{CD} and from Proposition \ref{TTT}.
\end{proof}
Thus, an $A_{\underline{k}}$-orbit of a rational map $f$ in a stratum $R_{\overline{m}}$ is biholomorphic to a $\prod_{i=1}^r A_{\underline{m^i}}$-orbit of an $\hat{f}\in \prod_{i=1}^r \text{\rm Rat}_{\underline{m^i}}\bigl(F(n)\bigr)$ such that $T(\hat{f})=f$. The same is true for the orbit structure of a fibre $\Pi^{-1}(c)$. Any $f\in \Pi^{-1}(c)$ belongs to the stratum $R_{\overline{m}}$, where the $\underline{m^i}=
(m^i_1,\dots,m^i_n)$ are determined by the polynomials \eqref{c_{ij}}, i.e.
\begin{equation} m^i_j=\text{order of vanishing of $q_j$ at $z_i$}.\label{order} \end{equation}
Combining this with \eqref{Pi_m} and \eqref{prod}, gives
\begin{proposition} Let $c\in \cx^{|{\underline{k}}|}$ and let $\{\underline{m^1},\dots,\underline{m^r}\}$ be the corresponding collection of multi-indices, determined by \eqref{order}. Then the $A_{\underline{k}}$-orbit structure of $\Pi^{-1}(c)$ is isomorphic to the $\prod_{i=1}^r A_{\underline{m^i}}$-orbit structure of $\Pi^{-1}(0,\dots 0)$, where this time $\Pi$ denotes the coordinate-wise projection \eqref{coord}.\hfill $\Box$\label{reduce}\end{proposition}
Thus, to study orbits of $A_{\underline{k}}$, we only need to know orbits of rational maps with a single pole located at $0$.

\section{Strongly regular points}
\par
By analogy with \cite{KW1}, we adopt the following definition:
\begin{definition} A point $f$ of $\Rat$ is said to be {\em strongly regular}, if its $A_{\underline{k}}$-orbit has maximal dimension ($=|\underline{k}|$), i.e. the isotropy group of $f$ is discrete.\label{sr}\end{definition}

According to Proposition \ref{reduce}, we only need to study orbits of rational maps with a single pole located at $0$. We shall now describe the structure of the set of strongly regular rational maps of this form.
\begin{theorem} Let $\underline{k}=(k_1,\dots,k_n)$ be a multidegree such that $k_i\geq 1$, $i=1,\dots,n$, and let $\Pi:\Rat\rightarrow \cx^{|\underline{k}|}$ be the pole location map \eqref{Pi}. The set of strongly regular points in $\Pi^{-1}(0)$ consists of $2^{n-1}$ orbits, each of which is isomorphic to
$\bigl(\cx^\ast\bigr)^n\times\prod_{i=1}^n \cx^{k_i-1}$.\label{N(c)}\end{theorem}

To remove the restriction $k_i\geq 1$, $i=1,\dots,n$, we observe that, if some $k_i=0$, then
$$ \Rat\simeq \text{\rm Rat}_{\underline{k^1}}\bigl(F(i-1)\bigr)\times \text{\rm Rat}_{\underline{k^2}}\bigl(F(n-i)\bigr),$$
where $\underline{k^1}=(k_1,\dots,k_{i-1})$, $\underline{k^2}=(k_{i+1},\dots k_n)$. Thus, the number of strongly regular orbits in an arbitraryl $\Pi^{-1}(0)$ is $2^t$, where $t=\#\{j;k_j\neq 0, k_{j+1}\neq 0\}$.

Combining this and Proposition \ref{reduce}, we obtain a description of the strongly regular part of any fibre $\Pi^{-1}(c)$:
\begin{corollary} Let $c\in \cx^{|{\underline{k}}|}$ and suppose that the associated polynomials $q_1(z),\dots,q_n(z)$, given by \eqref{c_{ij}}, have $r$ distinct zeros $z_1,\dots,z_r$. For every $i=1,\dots,r$, let 
$$s_i=\#\{j;\enskip q_j(z_i)=0\},\quad t_i=\#\{j;\enskip q_j(z_i)=0=q_{j+1}(z_i)\}$$ (with the convention that $q_{n+1}\equiv 1$) and let $s=\sum_{i=1}^r s_i$, $t=\sum_{i=1}^r t_i$. Then $\Pi^{-1}(c)$ contains 
$2^{t}$ $|\underline{k}|$-dimensional orbits of $A_{\underline{k}}$  and each such orbit is isomorphic to $\bigl(\cx^\ast\bigr)^s\times \cx^{|{\underline{k}}|-s}$.\hfill $\Box$\label{NN(c)}\end{corollary}

To prove Theorem \ref{N(c)}, we need first another description of $\Rat$.

\subsection{A matricial description of $\text{\rm Rat}_{\protect\underline{k}}$}
Let $\underline{k}=(k_1,\dots,k_n)$ be a multidegree. From the Appendix, we can infer a description of $\Rat$ as a symplectic quotient of $\prod_{i=1}^n T^\ast GL(k_i,\cx)$. This is given in Proposition \ref{BB} and Theorem \ref{Hurt}, and we describe now the end result. 
The key role is played by the following  generalised companion matices:
\begin{equation}B=\left(\begin{array}{ccc|cccc} &  &  & 0 &\ldots & 0 & b_1\\
& X & & \vdots & & \vdots & \vdots\\ &  &  & 0 &\ldots & 0 & b_m \\ \hline 
a_1 & \ldots & a_m & 0 &\ldots & 0 & c_1\\
0 &\ldots & 0 & 1 & \ddots & & c_2\\
\vdots & & \vdots & & \ddots & \ddots & \vdots \\0 &\ldots & 0 & 0&\ldots & 1 & c_{k-m} \end{array}\right).\label{B}\end{equation}
\begin{definition} Set $k_0=k_{n+1}=0$. We consider the following data:
to  every $i=1,\dots,n$, we associate a pair $B_i^-,B_i^+$ of $k_i\times k_i$ matrices and a $g_i\in GL(k_i,\cx)$  so that the following conditions are satisfied:
\begin{itemize}
\item $B_i^-$ (resp. $B_i^+$) is of the form \eqref{B} with $m=\min\{k_{i-1},k_i\}$ (resp. $m=\min\{k_{i+1},k_i\}$).
\item If $k_i>k_{i+1}$ (resp. if $k_i<k_{i+1}$), then the $X$-block of $B_i^+$ is equal to $B_{i+1}^-$ (resp. the $X$-block of $B_{i+1}^-$ is equal to $B_i^+$). 
\item If $k_i=k_{i+1}$, then there exists $(u_i,w_i)\in \cx^{k_i}\oplus \cx^{k_i} $ with $B_i^+-B_{i+1}^-=u_iw_i^T$.
\item $g_iB_i^+ g_i^{-1}=B_i^-.$
\end{itemize}
We denote by $\sM_{\underline{k}}$ the set of all such $F=\bigl(B_i^-,B_i^+,g_i,u_i,w_i\bigr)$. \label{sM}\end{definition}

The set $\sM_{\underline{k}}$ is acted upon by the group \begin{equation} G_{\underline{k}}=\prod_{i=1}^{n-1} GL(m_i,\cx),\enskip \text{where $m_i=\min\{k_{i},k_{i+1}\}$.}\label{group}\end{equation} To describe the action of an $h_i\in  GL(m_i,\cx)$, consider also such an $h_i$ as the following element of $GL(\max\{k_{i},k_{i+1}\},\cx)$:
\begin{equation} \left(\begin{array}{c|c} h_i & 0 \\ \hline 0 & 1\end{array}\right).\label{h_i}\end{equation}
Then 
\begin{equation} h_i:\enskip B_i^{+}\mapsto  h_iB_i^+h_i^{-1},\enskip B_{i+1}^{-}\mapsto  h_iB_{i+1}^- h_i^{-1}, \enskip g_i\mapsto g_ih_i^{-1},\enskip g_{i+1}\mapsto h_ig_{i+1}.\label{actor}\end{equation}
 Moreover, if $k_i=k_{i+1}$, then $h_i. (u_i,w_i^T)=\bigl(h_iu_i,w_i^Th_i^{-1}\bigr)$. The action is trivial on the remaining pieces of data.
The results described in  the Appendix can be now restated as follows:
\begin{proposition} $\Rat$ is biholomorphic to $\sM_{\underline{k}}/G_{\underline{k}}$. Moreover, the symplectic form of $\Rat$ is given by 
\begin{equation} \omega= \sum_{i=1}^n\tr\left(dg_ig_i^{-1}\wedge dB_i^-  -B_i^-(dg_ig_i^{-1}\wedge dg_ig_i^{-1})\right)-\sum_{\{j;k_j=k_{j+1}\}} dw_j^T\wedge du_j.\label{ooo}\end{equation}
The polar part $(q_1(z),\dots,q_n(z))$ of a rational map is $ q_i(z)=\det\bigl(z-B_i^-\bigr).$
\hfill $\Box$
\label{oooo}\end{proposition}

In particular, given Proposition \ref{LPP}, 
$$\LPP_{\underline{k}}\simeq \left\{\bigl(B_i^-,B_i^+,g_i,u_i,w_i\bigr)\in \sM_{\underline{k}}; \enskip \text{all $B_i^\pm$ are nilpotent}\right\}/G_{\underline{k}}.$$

\subsection{Proof of Theorem \ref{N(c)}}

We now consider the action of $A_{\underline{k}}$ in the above description. Let $p(z)=\sum_{j=1}^{k_i} \lambda_jz^j$ be a polynomial with null constant term, and consider the Hamiltonian $\tr p\bigl(B_i^-\bigr)$. The formula \eqref{ooo} implies that the corresponding Hamiltonian vector field is the right-invariant vector field induced by
$$ \frac{dp}{dz}\bigl(B_i^-\bigr)\in \gl(k_i,\cx),$$
and, hence, the element $(\lambda_1,\dots,\lambda_{k_i})\in \cx^{k_i}$ acts by 
\begin{equation} g_i\mapsto \exp\left( \frac{dp}{dz}\bigl(B_i^-\bigr)\right)g_i,\label{exp}\end{equation}
and it leaves unchanged every other piece of data. 
\par
We now represent $\lambda\in  A_{\underline{k}}\simeq\cx^{|\underline{k}|}$ by the coefficients of $n$ polynomials $p_i(z)=\sum_{j=1}^{k_i} \lambda_{ij}z^j$, so that the action of $\lambda$ is given by integrating the Hamiltonian vector field corresponding to $\sum_{i=1}^n\tr p_i\bigl(B_i^-\bigr)$.
Thus, the action of $\lambda$ on a rational map represented by $F=\bigl(B_i^-,B_i^+,g_i,u_i,w_i\bigr)\in \sM_{\underline{k}}$ is given by \eqref{exp}, and $\lambda$ belongs to the isotropy subgroup at $F$ if and only if there is an element $h=(h_1,\dots,h_n)$ of the group $G_{\underline{k}}$ given by \eqref{group} such that 
\begin{gather} h_iB_i^+h_i^{-1}=B_i^+,\enskip h_iB_{i+1}^- h_i^{-1}=B_{i+1}^-, \enskip h_iu_i=u_i,w_i^Th_i^{-1}=w_i^T,\label{gather1}\\ \exp\left( \frac{dp_i}{dz}\bigl(B_i^-\bigr)\right)g_i=h_{i-1}g_ih_i^{-1},\label{gather2}\end{gather}
for every $i=1,\dots,n$, where $h_0=1$. This gives {\em a} description of the orbits of $A_{\underline{k}}$. 
\par
We observe that \eqref{exp} implies:
\begin{lemma} A necessary condition for $\bigl(B_i^-,B_i^+,g_i,u_i,w_i\bigr)\in \sM$ to be strongly regular is that every $B_i^\pm$ is a regular matrix.\hfill $\Box$\end{lemma}
To prove Theorem \ref{N(c)}, we only need to consider rational maps with a single pole located at $0$. Because of this and the last lemma,  we can assume that all $B_i^\pm$ are regular nilpotent matrices. Any such a matrix is conjugate to
\begin{equation}\left(\begin{array} {ccccl}0 & 0 &\dots & 0 & 0\\ 1 & 0 &\dots & 0 & 0\\0 & 1 &\dots & 0 & 0\\
\vdots &\vdots & \ddots &\vdots &\vdots\\0 & 0 &\dots & 1& 0\end{array}\right).\label{nilp}\end{equation}
We now use the action of $G_{\underline{k}}$, given by \eqref{group}, to make $B_i^\pm$ of a the following form:
\begin{itemize} 
\item If $k_i>k_{i+1}$, then $B_{i+1}^-$ is of the form \eqref{nilp}.
\item If $k_i\leq k_{i+1}$, then $B_{i}^+$ is of the form \eqref{nilp}.
\end{itemize}
Definition \ref{sM} imposes now conditions on the remaining $B_i^\pm$. For $k_i\neq k_{i+1}$ the larger of the two matrices $B_i^+,B_{i+1}^-$ is of the form \eqref{B} with $X$ being \eqref{nilp}. Since, in the formula \eqref{B}, we have
$$ \det(z-B)=\det(z-X)\bigl(z^{k-m}-c_{k-m}z^{k-m-1}+\dots+(-1)^{k-m}c_1\bigr)+a(z-X)_{\rm adj}b,$$
where the subscript {\em adj} denotes the classical adjoint matrix, we get, given that $\det(z-B)=z^k$ and $\det(z-X)=z^m$, that 
$e_1=\dots=e_{k-m}=0$ and $f(z-X)_{\rm adj}g=0$. The second condition implies, in particular, that
\begin{equation} a_mb_1=0.\label{fg}\end{equation}
Similarly, when $k_i=k_{i+1}$, then $B_{i+1}^-=B_i^+-uw^T$, and, since $B_i^+$ is of the form \eqref{nilp}, we obtain from the formula 
$$ \det\bigl(z-B_{i+1}^-\bigr)=\det\bigl(z-B_i^+\bigr)-w^T\bigl(z-B_i^+\bigr)_{\rm adj} u$$
that $w^T(z-X)_{\rm adj} u=0$, where $X$ is the matrix \eqref{nilp}. Again, we have
\begin{equation} w_{k_i}u_1=0.\label{uw}\end{equation}
To a rational map, given by $F=\bigl(B_i^-,B_i^+,g_i,u_i,w_i\bigr)\in \sM_{\underline{k}}$, we now associate a map $\sigma_F:\{1,\dots,n-1\}\rightarrow \{-1,0,1\}$ by setting, for every $i=1,\dots n-1$, 
$$ \sigma_F(i)=\begin{cases} -1 & \text{if $k_i\neq k_{i+1}$ and $a_m\neq 0$ in  \eqref{fg} or if $k_i= k_{i+1}$ and $w_{k_i}\neq 0$ in  \eqref{uw}}\\
0 & \text{ if $k_i\neq k_{i+1}$ and $a_m=b_1= 0$ or $k_i= k_{i+1}$ and $w_{k_i}=u_1= 0$}\\
1 & \text{if $k_i\neq k_{i+1}$ and $b_1\neq 0$ in  \eqref{fg}  or if $k_i= k_{i+1}$ and $u_1\neq 0$ in  \eqref{uw}}.\end{cases}$$
Theorem \ref{N(c)} follows from:
\begin{proposition} Suppose that all degrees $k_i$, $i=1,\dots,n$, are nonzero. A rational map $F$ in $\Pi^{-1}(0)$ is strongly regular if and only if the map $\sigma_F$ does not take value $0$. If $F$ is strongly regular, then its $A_{\underline{k}}$-orbit is of the form stated in Theorem \ref{N(c)}. Moreover, assigning $F\mapsto \sigma_F$ sets  a $1-1$ correspondence between orbits of strongly regular points in $\Pi^{-1}(0)$ and maps $\sigma:\{1,\dots,n-1\}\rightarrow \{-1,0,1\}$ which do not take value $0$.\label{sigma}\end{proposition} 
\begin{proof} Observe that, since we assume that, for every $i=1,\dots,n-1$, either $B_i^+$ or $B_{i+1}^-$ is of the form \eqref{nilp}, we can conclude that there exists a non-trivial $h_i\in GL(\min\{k_i,k_{i+1}\},\cx)$ satisfying \eqref{gather1} if and only if $\sigma_F(i)=0$. Thus, if  $\sigma_F$ does not take value $0$, then $F$ is strongly regular, and the formula \eqref{exp} implies that its orbit is of the stated form. We now show that, if $\sigma_F(i)=0$, then $F$ is not strongly regular. Since $\sigma_F(i)=0$, any element $h_i\in GL(\min\{k_i,k_{i+1}\},\cx)$ of the form
$$\left(\begin{array} {cccc}1 & 0 &\dots & 0 \\ 0 & 1 &\dots & 0\\
\vdots &\vdots & \ddots &\vdots \\n & 0 &\dots & 1\end{array}\right)$$
satisfies \eqref{gather1}. In particular,
$h_i$ commutes with both $B_i^+$ and $B_{i+1}^-$ ($h_i$ is viewed as an element of $GL(\max\{k_i,k_{i+1}\},\cx)$ via \eqref{h_i}). Since $B_i^+$ and $B_{i+1}^-$ are regular nilpotent matrices, we can find polynomials $p_i(z),p_{i+1}(z)$, of degrees $k_i$ and $k_{i+1}$, such that 
$$h_i=\exp\left( \frac{dp_i}{dz}\bigl(B_i^+\bigr)\right),\enskip h_i=\exp\left( \frac{dp_{i+1}}{dz}\bigl(B_{i+1}^-\bigr)\right).$$
Then the element $(1,\dots,1,h_i,1,\dots,1)$ of $G$ satisfies \eqref{gather1} and \eqref{gather2}, and we get a continuous subgroup of the isotropy group. $F$ cannot be strongly regular.
\par
We now show that for any choice of a map $\sigma:\{1,\dots,n-1\}\rightarrow \{-1,1\}$, all $F$ with $\sigma_F=\sigma$ belong to a single $A_{\underline{k}}$-orbit. We observe that the action of $A_{\underline{k}}$ leaves invariant the map $\sigma_F$. This follows from  \eqref{exp} and our chosen form of each $B_i^-$. Therefore, rational maps with different $\sigma_F$ belong to different $A_{\underline{k}}$-orbits. We also observe, that, for $i=1,\dots,n-1$, if $a_m\neq 0$ (resp. $b_1\neq 0$), then there is a unique $h_i$ centralising the matrix \ref{nilp}, such that $ah_i^{-1}=(0,\dots,0,1)$ (resp. $h_ib=(1,0,\dots,0)^T$). Similarly, when $k_i=k_{i+1}$, we can make $w=(0,\dots,0,1)^T$ or $u=(1,0,\dots,0)^T$. This gets rid of the remaining action of $G_{\underline{k}}$. Moreover, any $B_i^\pm$ is now either the matrix \eqref{nilp}
(when $a=(0,\dots,0,1)$ or $w=(0,\dots,0,1)$) or the matrix $E$ which differes from \eqref{nilp} by having $1$ as the $(1,k_i)$-entry and $0$ as the $(m+1,m)$-entry ($m$ refers to the form \eqref{B}). Therefore, any two $g_i$, which conjugate $B_i^+$ to $B_i^-$, differ by a polynomial in $B_i^-$ (as these are regular nilpotent). Formula \eqref{exp} shows now that $F,F^\prime$ with the same $\sigma$ belong to the same orbit of $A_{\underline{k}}$.
\par
It remains to prove that any $\sigma:\{1,\dots,n-1\}\rightarrow \{-1,1\}$ is realised as $\sigma_F$ for a strongly regular $F\in \Pi^{-1}(0)$. In the previous paragraph, we showed that for any $\sigma:\{1,\dots,n-1\}\rightarrow \{-1,1\}$, there is a canonical form of the $B_i^\pm$, which are either \eqref{nilp} or the matrix $E$. Since $E$ is conjugate to \eqref{nilp}, we get elements $g_i\in GL(k_i,\cx)$ such that $g_iB_i^+g_i^{-1}=B_i^-$, and so, according to Proposition \ref{oooo} and Definition \ref{sM}, a rational map $F$ with $\sigma_F=\sigma$.
\end{proof}

\part{Gelfand-Zeitlin actions}

\section{Gelfand-Zeitlin actions on manifolds\label{T*GL}}

As mentioned in the introduction, an analogue of a Gelfand-Zeitlin action exists on any symplectic manifold with an Hamiltonian $GL(n,\cx)$-action. We shall explain this now, beginning with the symplectic structure of $T^\ast GL(n,\cx)$.

\subsection{The symplectic structure of  $T^\ast G$} Let $G$ be a real or complex semisimple Lie group with Lie algebra $\fG$. Its cotangent bundle $T^\ast G$ can be trivialised using right-invariant differential forms. We can further identify $\g^\ast\simeq \g$ using the Killing form and so we obtain an isomorphism $T^\ast G\simeq G\times \fG$.
We consider the canonical symplectic form $\omega$ of $T^\ast G$ in this trivialisation. At a point $(g,B)\in G\times \fG$, $\omega$ is given by
\begin{equation} \omega\bigl((\rho,b),(\hat{\rho},\hat{b})\bigr)=\tr \bigl(\rho \hat{b}-\hat{\rho}b-B[\rho,\hat{\rho}]\bigr),\label{omega}\end{equation}
where $\rho,\hat{\rho}$ are right-invariant vector fields on $G$, $b,\hat{b}$ are translation-invariant vector fields on $\g$ and $\tr XY$ denotes Killing form on $\g$.
\par
The symplectic structure is invariant under both left and right action by $G$ on $T^\ast G$. In the chosen trivialisation $T^\ast G\simeq G\times \fG$, the actions are given by
\begin{equation}h \cdot_L (g,B)=\bigl(hg,\Ad(h)B\bigr),\quad h \cdot_R (g,B)=\bigl(gh^{-1},B\bigr),\label{action}\end{equation}
and the corresponding moment maps are
\begin{equation} \mu_L(g,B)=B ,\quad \mu_R(g,B)=-\Ad\bigl(g^{-1}\bigr)B.\label{moment}\end{equation}

\begin{remark} In the category of Hamiltonian  $G$-spaces, $T^\ast G$ acts as an identity for the symplectic quotient, i.e. if $M$ is an Hamiltonian $G$-space, then the symplectic quotient of $M\times T^\ast G$ by the diagonal $G$ (with the right action on $T^\ast G$) is canonically isomorphic to $M$.\label{triv}\end{remark}
\par
We need to consider one more subgroup of  symplectomorphisms of $T^\ast G$: the group $S$ obtained by integrating  Hamiltonian vector fields corresponding to $G\times G$-invariant functions. Any $G\times G$-invariant function on $T^\ast G$ is given by a unique $G$-invariant function on $\g$ only, and can be identified with a $W$-invariant function on a Cartan subalgebra $\fH$.
\par
For $G=GL(m,\cx)$ we shall denote by $C_m$ the subgroup of $S$ generated by the Hamiltonians  $\gl(m,\cx)\ni B\mapsto \tr B^i$, $i=1,\dots, m$.  The action of $C_m$ is made more explicit by
\begin{proposition} The Hamiltonian vector field $X_j$, generated by $\gl(m,\cx)\ni B\mapsto \tr B^j$, on $T^\ast GL(m,\cx)\simeq GL(m,\cx)\times \gl(m,\cx)$  is equal, at the point $(g,B)$,  to the right-invariant vector field generated by $j B^{j-1}$, and, hence, it induces a global $\cx$-action, given by 
$$ z.(g,B)=\left(e^{jzB^{j-1}}g,B\right)$$.\label{Ham}\end{proposition}
\begin{proof} This is a simple calculation using \eqref{omega}.\end{proof}

Thus, $C_m\simeq \cx^m$ is abelian  and its action commutes with left and right actions of $GL(m,\cx)$ on $T^\ast GL(m,\cx)$.

\subsection{The Gelfand-Zeitlin action on $GL(n,\cx)$-manifolds\label{man}}

All manifolds, discussed here, are complex and all symplectic structures, group actions, etc. are holomorphic. All symplectic quotients are taken at the zero-level of relevant moment maps.

\begin{proposition} Let $M$ be a symplectic manifold with an Hamiltonian $GL(n,\cx)$-action. Then $M$ is isomorphic, as a symplectic manifold, to the symplectic quotient of 
\begin{equation}M\times  T^\ast GL(n,\cx)\times\cdots\times T^\ast GL(1,\cx)\quad \text{by}\enskip H=\prod_{i=1}^{n} GL(i,\cx),\label{quotient}\end{equation}
where each $GL(i,\cx)$, $i<n$, acts on the right of  $T^\ast GL(i,\cx)$ and on the left on $T^\ast GL(i+1,\cx)$ as the group of upper-left $i\times i$-minors. The group $GL(n,\cx)$ acts in a given way on $M$ and on the right of $T^\ast GL(n,\cx)$.\label{trivial}
\end{proposition}
\begin{proof} This follows by induction from Remark \ref{triv}.\end{proof}
\begin{remark} It also follows from  Remark \ref{triv} that $T^\ast GL(n,\cx)$  and the group $GL(n,\cx)$ can be omitted in \eqref{quotient}.\label{omit}\end{remark}

We now observe that on $\prod_{i=1}^{n} T^\ast GL(i,\cx)$ we have the Hamiltonian action of abelian group $A=\prod_{i=1}^n C_i$ and, since $A$ commutes with $\prod_{i=1}^{n} GL(i,\cx)$, we obtain a Hamiltonian action of $A$ on $M$. We identify this action as follows:
\begin{proposition} Let $a_{mi}$ be the element of $C_m$ corresponding to  $\gl(m,\cx)\ni B\mapsto \tr B^i$. Then the Hamiltonian of $a_{mi}$ on $M$ is given (up to an additive constant) by 
$$ \mu_{mi}(p)=\tr \left( \mu(p)_{(m)}\right)^i,$$
where $\mu:M\rightarrow \gl(n,\cx)$ is the moment map and $B_{(m)}$ denotes the upper-left $m\times m$-minor of a matrix $B$.\label{minor}\end{proposition}
\begin{proof}  We view $M$ as in \eqref{quotient}, i.e. as a symplectic quotient of $M\times \prod_{j=1}^{n} T^\ast GL(j,\cx)$. We write elements of $T^\ast GL(j,\cx)$ as $(g_j,B_j)$ in the right trivialisation.  The moment map equations  are
$$ \mu(p)-g_n^{-1}B_ng_n=0,\enskip (B_{j+1})_{(j)}-g_j^{-1}B_jg_{j}=0\enskip \forall j\leq n-1.$$
Taking the quotient corresponds now to acting by the element $\bigl(g_n^{-1},\dots,g_1^{-1}\bigr)$. Thus, the quotient is identified with the set of points $\bigl(p,(B_n,1),\dots,(B_1,1)\bigr)\in M\times \prod_{j=1}^{n} T^\ast GL(j,\cx)$ satisfying:
$$ \mu(p)=B_n, \enskip \enskip (B_{j+1})_{(j)}=B_j\enskip \forall j\leq n-1.$$ 
Hence, $\tr \bigl(\mu(p)_{(m)}\bigr)^i=\tr (B_m)^i$ and the latter is the Hamiltonian for $a_{mi}$ acting on $T^\ast GL(m,\cx)$. Since this Hamiltonian is invariant under the group $H$, it descends to the symplectic quotient, which proves the statement.\end{proof}

We conclude 
\begin{corollary}Let $M$ be a (complex) symplectic manifold with an Hamiltonian $GL(n,\cx)$-action, and let $\mu:M\rightarrow \gl(n,\cx)$ be the moment map. For a point $p\in M$, denote by $B_{mi}(p)$ the element $\left( \mu(p)_{(m)}\right)^{i-1}$ of the algebra $\gl(n,\cx)$, $m\leq n$, $i\leq m$. Then the actions of $\cx $ on $M$
$$z. p=e^{zB_{mi}(p)}. p$$
commute for various $m$ and $i$ and induce an action of $\cx^{n(n+1)/2}$ on $M$. \label{GZ1}\end{corollary}
For obvious reasons we call this action on $M$ the {\em Gelfand-Zeitlin action}. 
\par
We can, in particular, take $M=T^\ast GL(n,\cx)$ with the left action. 
 The last proposition and the commutativity of the action of $A$ and $GL(n,\cx)$ give immediately
\begin{corollary} The action of $A$ on $T^\ast GL(n,\cx)$ descends to the Poisson quotient $T^\ast GL(n,\cx)$ by the right $GL(n,\cx)$, i.e. to $\gl(n,\cx)^\ast$ and coincides there with the usual Gelfand-Zeitlin action.\hfill $\Box$\label{matrices}\end{corollary}

\begin{remark} Corollary \ref{GZ1} remains valid  for any holomorphic Poisson manifold $M$ equipped with a complete Poisson map $\mu:M\rightarrow \gl(n,\cx)^\ast$. Indeed, the Gelfand-Zeitlin action is defined in the same way, and the commutativity of various Hamiltonians follows from applying Corollary \ref{GZ1} to the symplectic leaves of $M$. Alternatively, once we have the commutativity for the action on $\gl(n,\cx)^\ast$, the commutativity on $M$ follows.\end{remark}

\section{The symplectic manifold $GL(n,\cx)\times \cx^n$\label{Vn}}

\subsection{The symplectic structure}
We consider now the manifold defined by \eqref{V_n}, i.e. 
\begin{equation} V_n=\{(B,b)\in \gl(n,\cx)\times\cx^n; \enskip\text{$b$ is a cyclic vector for B}\}.\label{Bb}\end{equation}
It is the set of stable points of the  action $GL(n,\cx)$ on $\gl(n,\cx)\times\cx^n$ given by
\begin{equation}g.(B,b)=\bigl(gBg^{-1},gb\bigr).\label{action2}\end{equation}
$V_n$ is biholomorphic to $GL(n,\cx)\times \cx^n$ via the map
\begin{equation} (B,b)\mapsto \left(\bigl(b,Bb,\dots,B^{n-1}b\bigr)\,,\,\bigl(\tr B,\dots,\tr B^{n}\bigr)\right), \label{iso}\end{equation}
and it is in this guise that the symplectic structure of $V_n$ is easiest seen. One has the obvious isomorphism between  $GL(n,\cx)\times \cx^n$ and $GL(n,\cx)\times S$, where $S$ is the set of companion matrices, i.e. matrices of the form
\begin{equation} \left(\begin{array}{ccccl}0 & 0 &\dots & 0 & \beta_n\\ 1 & 0 &\dots & 0 & \beta_{n-1}\\0 & 1 &\dots & 0 & \beta_{n-2}\\
\vdots &\vdots & \ddots &\vdots &\vdots\\0 & 0 &\dots & 1& \beta_1\end{array}\right).\label{comp}\end{equation}
Now, the set $GL(n,\cx)\times S $ can be realised (see, e.g., \cite{BielLMS}) as the symplectic quotient of $T^\ast GL(n,\cx)$ by the left action \eqref{action} of the group $N^+$ of upper-triangular unipotent matrices. The relevant level of the moment map is the following invariant element of $\Lie(N^+)^\ast\simeq \Lie(N^-)$:
$$ \left(\begin{array}{ccccl}0 & 0 &\dots & 0 & 0\\ 1 & 0 &\dots & 0 & 0\\0 & 1 &\dots & 0 & 0\\
\vdots &\vdots & \ddots &\vdots &\vdots\\0 & 0 &\dots & 1& 0\end{array}\right).$$
The symplectic form of $V_n$ is therefore given, at a point $(g,C)\in GL(n,\cx)\times S$, by the formula \eqref{omega}. It is invariant under the right action of $GL(n,\cx)$ on itself, and the moment map, from \eqref{moment}, is:
$$\mu(g,C)=-\Ad\bigl(g^{-1}\bigr)C.$$
In the original description \eqref{Bb} of $V_n$, the moment map is
\begin{equation}\mu(B,b)=B.\label{moma}\end{equation}

We finish this section by remarking that $V_n$ has a geometric interpretation \cite{Craw} as the space of {\em full} rational maps $\oP^1\rightarrow \oP^n$ of degree $n$. Here ``full" means that the image of the rational map is not contained in any hyperplane. We shall not, however, make any use of this fact.

\subsection{The Gelfand-Zeitlin action}  We can apply Corollary \ref{GZ1} to $V_n$, with its $GL(n,\cx)$ action and the moment map \eqref{moma}, to obtain:

\begin{proposition} The $n(n+1)/2$ functions
 $$\mu_{mi}:V_n\rightarrow \cx,\quad \mu_{mi}(B,b)=\tr ( B_{(m)})^i,$$
where $B_{(m)}$ denotes the upper-left $m\times m$-minor of $B$,
are Poisson-commuting on $V_n$ and their Hamiltonian vector fields are globally integrable defining a symplectic action of a group $A\simeq \cx^{n(n+1)/2}$ on $V_n$.\hfill $\Box$\end{proposition} 
We observe that the subgroup $C_n\simeq \cx^n$ of $A$, given by integrating the flows corresponding to $\mu_{n1},\dots,\mu_{nn}$, has, as the moment map, the characteristic polynomial of $B$, and therefore the symplectic quotients of $V_n$ by $C_n$ are the regular adjoint orbits of $GL(n,\cx)$.
\par
The action of $A$ has the following geometric interpretation, already stated as Theorem \ref{theorema} of the introduction. Recall the definition of the group  $A_{\underline{k}}$ and the maps $c_{ij}$ from Section \ref{A(c)}.
\begin{theorem} There exists a canonical symplectomorphism $\Phi:V_n\rightarrow \sR_n$, where  $\sR_n=\text{\rm Rat}_{(1,\dots,n)}\bigl(F(n)\bigr)$ is the space
of based rational maps of degree $(1,\dots,n)$ from $\oP^1$ to the manifold of full flags in $\cx^{n+1}$. Moreover, $\Phi$ is equivariant for the actions of $A$ on $V_n$ and $A_{(1,\dots,n)}$ on $\sR_n$ and it intertwines the relevant moment maps, i.e. $\Phi^\ast c_{mi}=\mu_{mi}$ for every $m\leq n$ and $i\leq m$.  \label{Phi}\end{theorem}
\begin{proof} The discussion in Section \ref{man} shows that the Gelfand-Zeitlin action on $V_n$ arises from identification of $V_n$ with the symplectic quotient of $V_n\times   T^\ast GL(n-1,\cx)\times\cdots\times T^\ast GL(1,\cx)$ by $\prod_{i=1}^{n-1} GL(i,\cx)$ (as observed in Remark \ref{omit}, the factor $T^\ast GL(n,\cx)$ can be omitted). This (complex) symplectic quotient coincides, given \eqref{isomorphism1} and \eqref{isomorphism2}, with the symplectic quotient described in Corollary \ref{Phi1} of the Appendix. Thus, we obtain a symplectomorphism $\Phi$ between $V_n$ and $\sR_n$. Since the groups $A$ and $A_{(1,\dots,n)}$ are defined by integrating Hamiltonian vector fields corresponding to $c_{mi}$ and $\mu_{mi}$, it is enough to show that $\Phi^\ast c_{mi}=\mu_{mi}$ for every $(m,i)$. This follows, however, from the above discussion and the formula \eqref{plar}.
\end{proof}

We are going to give two applications of this theorem. For the first one, we recall that Kostant and Wallach introduced in \cite{KW1} the concept of {\em strongly regular matrices}, i.e. matrices whose stabiliser for the Gelfand-Zeitlin action is discrete.  Theorem 3.14 in \cite{KW1} describes the Gelfand-Zeitlin orbit of a strongly regular matrix and agrees with our Corollary \ref{NN(c)}. However, we can now strengthen results of \cite{KW1} to determine  the number of  strongly regular orbits with prescribed  spectral data:
\begin{corollary} Let $c=(c_{mi})_{i\leq m}\in \cx^{n(n+1)}$ and let $M_c(n)$ be the set of all $n\times n$ complex matrices $B$ such that the characteristic polynomial of the $m\times m$ upper-left minor of $B$ has $c_{m1},\dots,c_{mm}$ as coefficients, i.e.
\begin{equation} \chi_m:=\det\bigl(z-B_{(m)}\bigr)=z^m+\sum_{i=1}^{m} c_{mi}z^{i-1}, \enskip m=1,\dots,n.\label{chi}\end{equation}
Suppose that the polynomials $\chi_1,\dots,\chi_{n}$ have $r$ distinct roots $z_1,\dots, z_r$ and set $t_i=\#\{j;\enskip \chi_j(z_i)=0=\chi_{j+1}(z_i) \}$. Finally,  let $t=\sum_{i=1}^r t_i$. Then  $M_c(n)$ contains exactly $2^{t}$ orbits of maximal dimension. 
\end{corollary}
\begin{proof} First of all, since we are interested only in strongly regular, hence regular matrices, we can restrict ourselves to the subset $M_c^\text{\rm reg}(n)\subset M_c(n)$ consisting of regular matrices. This subset is the quotient of $\Pi^{-1}(c)\subset \text{Rat}_{(1,\dots,n)}\bigl(F(n)\bigr)$ by the subgroup generated by the Hamiltonians $c_{n1},\dots,c_{nn}$. More generally, for any $\underline{k}=(k_1,\dots,k_n)$, consider the quotients $M_d^\text{\rm reg}(\underline{k})$ of $\Pi^{-1}(d)\subset \Rat$
by the subgroup generated by the Hamiltonians $d_{n1},\dots,d_{nk_n}$.  Now, observe that the decomposition of the action, given in Proposition \ref{reduce} commutes with taking these quotients. Therefore, we need only to know the number of orbits of maximal dimension in $M_0^\text{\rm reg}(\underline{k})$, for any $\underline{k}$. Proposition \ref{sigma} implies that this number is the same as the number of orbits of maximal dimension in $\Pi^{-1}(0)$. Thus, the number of orbits of maximal dimension in $M_c(n)$ is the one given in Corollary \ref{NN(c)}.
\end{proof}

\medskip

The second application of Theorem \ref{Phi} is the observation, that the constructions and results of \cite{KW2} have a very natural interpretation, given the isomorphism $\Phi$.  Kostant and Wallach consider the subset $M_\Omega(n)$ of $n\times n$ matrices such that the polynomials \eqref{chi} have all roots distinct, i.e. all the eigenvalues of all main minors are distinct. They consider then a covering $M_\Omega(n,\fE)$ of $M_\Omega(n)$, whose group of deck transformations is the direct product of symmetric groups $S_m$, $m=1,\dots,n$. They show  that   the roots $r_l$ of the polynomials \eqref{chi} form a maximal commutative subalgebra of  $\sO\bigl(M_\Omega(n,\fE)\bigr)$ for the lifted Poisson structure. Kostant and Wallach construct then dual coordinates $s_l$ (also commuting) and show that 
\begin{equation} [r_l,s_k]=\delta_{lk}s_{k}.\label{commutor}\end{equation}
\par
We find that this picture becomes very clear under the diffeomorphism $\Phi$ of Theorem \ref{Phi}. The subset  $M_\Omega(n)$ coincides now with the subset $p_n\equiv 1$ of the open stratum $\Rat^\circ$ described at the beginning of section \ref{symplectic} and $M_\Omega(n,\fE)$ is the set of $\bigl\{z_j^i,p_i\bigl(z_j^i\bigr)\bigr\}$ given in \eqref{sympl}, i.e. the subset of $\prod \text{\rm Rat}_1\bigl(\oP^1\bigr)$ where all poles are distinct (with $p_n(z^n_j)=1$ for $j=1,\dots,n$).  The Kostant-Wallach coordinates $r_l$ are the $z_j^i$ and the $s_l$ are $p_i\bigl(z_j^i\bigr)^{-1}$. The relation \eqref{commutor} is now easily seen from \eqref{sympl}.

\section{Gelfand-Zeitlin actions on $T^\ast GL(n,\cx)$}

As discussed in section \ref{T*GL},  $T^\ast GL(n,\cx)$ is a symplectic manifold with two, left and right, actions of $GL(n,\cx)$. Since these actions commute, Corollary \ref{GZ1} gives an action of two commuting copies of $\cx^{n(n+1)/2}$ on $T^\ast GL(n,\cx)$. These are induced by integrating Hamiltonian vector fields corresponding to upper-left minors of the moment maps for the left and right actions, i.e., according to \eqref{moment}, to the minors of $B$ and $-g^{-1}Bg$, where $(g,B)\in GL(n,\cx)\times \gl(n,\cx)\simeq T^\ast GL(n,\cx)$.  In particular, the actions induced by the characteristic polynomials of the full matrices $B$ and $-g^{-1}Bg$ coincide, and, hence, the group acting effectively  on $T^\ast GL(n,\cx)$ is  $\tilde{A}\simeq \cx^{n^2}$.
\begin{remark} The Hamiltonians for $\tilde{A}$ are considered in \cite{GKL}, in the case of $T^\ast GL(n,\oR)$, as a maximal commutative subalgebra of the  Poisson algebra of functions on $T^\ast GL(n,\oR)$.\end{remark}

We can now identify the action of $\tilde{A}$ on  $T^\ast GL(n,\cx)$ with an action of the group $A_{\underline{k}}$ (see section \ref{A(c)}) on the appropriate space of rational maps:
\begin{theorem} There exists a canonical symplectomorphism $\Psi:T^\ast GL(n,\cx)\rightarrow \tilde{\sR}_{2n-1}$, where  $\tilde{\sR}_{2n-1}=\text{\rm Rat}_{(1,\dots,n,\dots,1)}\bigl(F(2n-1)\bigr)$ is the space
of based rational maps of degree $(1,2,\dots,n-1,n,n-1,\dots,2,1)$ from $\oP^1$ to the manifold of full flags in $\cx^{2n}$. Moreover, $\Psi$ is equivariant for the actions of $\tilde{A}$ on $T^\ast GL(n,\cx)$ and $A_{\underline{k}}$, $\underline{k}=(1,2,\dots,n-1,n,n-1,\dots,2,1)$, on $\tilde{\sR}_{2n-1}$ and it intertwines the relevant moment maps, i.e. $\Psi^\ast c_{mi}=\mu_{mi}$ for every $m\leq 2n-1$ and $i\leq\min\{m,2n-m\}$.  \label{Psi}\end{theorem}
\begin{proof} It is the same proof as for Theorem \ref{Phi}, using Corollary \ref{Psi1} instead of  Corollary \ref{Phi1}.\end{proof}

While the  group $\tilde{A}$ is a maximal commutative subgroup of the group of algebraic symplectomorphisms of $T^\ast GL(n,\cx)$, Corollary \ref{matrices} implies that for some purposes it would be better to describe the action of $A$ corresponding only to the left action of $GL(n,\cx)$ on  $T^\ast GL(n,\cx)$. It turns out, as we now explain, that there is a symplectomorphism, equivariant with respect to $A\times GL(n,\cx)$, between $T^\ast GL(n,\cx)$ and a space of rational maps into the manifold of {\em half-full flags} in $\cx^{2n}$.
\par
To see this, informally at first, consider  $T^\ast GL(n,\cx)$ as a symplectic quotient of $ T^\ast GL(n,\cx)\times T^\ast GL(n-1,\cx)\times\dots \times T^\ast GL(1,\cx)$ by $\prod_{i=1}^{n-1} GL(n-1,\cx)$. According to the discussion in the Appendix, in particular \eqref{isomorphism1} and \eqref{isomorphism2}, this space is equivariantly isomorphic to the moduli space of solutions to the Lax equation \ref{Lax} on the union of $n$ adjoining intervals, with  solutions $\bigl(\alpha^i(t),\beta^i(t)\bigr)$ on the $i$-th interval being $\gl(i,\cx)$-valued and such that, for $i<n$,  $\bigl(\alpha^{i}(t_{i+1}),\beta^i(t_{i+1})\bigr)$ coincides with  the upper-left $i\times i$ minor of $\bigl(\alpha^{i+1}(t_{i+1}),\beta^{i+1}(t_{i+1})\bigr)$. There are no poles at $t_{n+1}$ this time. Comparing this with Proposition \ref{Psi1}, we see that this corresponds to collapsing the union of the remaining intervals $[t_{n+i},t_{n+i+1}]$ to the point $t_{n+1}$. According to  Remark \ref{t}, this corresponds to replacing the adjoint $U(2n)$-orbit of $\sqrt{-1}(t_1,\dots,t_{2n})$ with the orbit of $\sqrt{-1}(t_1,\dots,t_{n+1},t_{n+1},\dots,t_{n+1})$, i.e. with the manifold $GL(2n,\cx)/P$ of flags of type
$E_1\subset E_2\subset\dots\subset E_n\subset \cx^{2n}$ with $\dim E_i=i$. Thus, $T^\ast GL(n,\cx)$ should be $A\times GL(n,\cx)$-equivariantly isomorphic to a subset of based rational maps to $GL(2n,\cx)/P$ of degree $(1,\dots,n)$. We shall now prove this and identify the relevant subset.
\par
First of all, we now think, as in \cite{Hu}, of rational maps as flags of subundles of a trivial bundle on $\oP^1$. Thus, a based rational map of degree $(k_1,\dots,k_m)$ into $F(m)$ corresponds to a flag $E_1\subset E_2\subset\dots\subset E_n\subset \sO_{\oP^1}\otimes \cx^{m+1}$ of subbundles such that the rank of $E_i$ is $i$, $E_i/E_{i-1}\simeq \sO(k_{i-1}-k_i)$, and $(E_i)$ coincides with the anti-standard flag at infinity.
\par
Similarly, a based rational map of degree $(l_1,\dots,l_n)$ into the manifold $GL(2n,\cx)/P$ of half-full flags corresponds to a flag $E_1\subset \dots\subset E_n\subset \sO^{2n}$ of subbundles such that the rank of $E_i$ is $i$, $E_i/E_{i-1}\simeq \sO(l_{i-1}-l_i)$, and $(E_i)$ coincides with the first half of the anti-standard flag at infinity. This space of rational maps into $GL(2n,\cx)/P$
has a natural stratification (cf. \cite{Murr}) according to the splitting of the rank $n$ vector bundle $\sO^{2n}/E_n$:
$$\sO^{2n}/E_n\simeq \bigoplus_{i=1}^n \sO(m_i),\quad \sum_{i=1}^n m_i=l_n.$$
Let us denote this stratum by $\text{\rm Rat}_{\underline{l};\underline{m}}\bigl(GL(2n,\cx)/P\bigr)$. We observe that the natural projection
$$\pi:F(2n-1)\rightarrow GL(2n,\cx)/P$$
induces a holomorphic map
\begin{equation} \pi_\ast: \text{\rm Rat}_{\underline{k}}\bigl(F(2n-1)\bigr)\rightarrow \text{\rm Rat}_{\underline{l};\underline{m}}\bigl(GL(2n,\cx)/P\bigr),\label{pi_star}\end{equation}
where $\underline{l}=(k_1,\dots,k_n)$, $\underline{m}=(k_n-k_{n+1},\dots,k_{2n-2}-k_{2n-1},k_{2n-1}).$
In general, this map is not surjective. We have, however:
\begin{proposition} Let $(k_1,\dots,k_{2n-1})$ satisfy
 $$k_n-k_{n+1}=k_{n+1}-k_{n+2}\dots=k_{2n-2}-k_{2n-1}=k_{2n-1}.$$
Then the map \eqref{pi_star} is a biholomorphism.\label{SSSS}\end{proposition}
\begin{proof} We construct an inverse of $\pi_\ast$. Let $f\in \text{\rm Rat}_{\underline{l};\underline{m}}\bigl(GL(2n,\cx)/P\bigr)$ correspond to a flag $E_1\subset \dots\subset E_n\subset \sO^{2n}$ as above. We wish to extend this flag to a full flag of subbundles, coinciding with the anti-standard flag at infinity.
\par
By assumption, $\sO^{2n}/E_n\simeq \sO(m)^n$, where $m=k_n-k_{n+1}=\dots= k_{2n-1}-0$. Moreover, at infinity, $\sO^{2n}/E_n$ is identified with the subspace of $\cx^{2n}$ spanned by $e_1,\dots,e_n$ (see the beginning of Section \ref{rff}). The bundle $\sO(m)^n$ also has a canonical trivialisation at infinity. Since $GL(n,\cx)$ acts on the set of isomorphisms $\phi:\sO^{2n}/E_n\simeq \sO(m)^n$, there is a unique $\phi$, which is the identity at infinity. We now set $E_{n+i}=p^{-1}\phi^{-1}\bigl(\sO(m)^i\bigr)$, where $\sO(m)^i\subset \sO(m)^n$ is spanned by the last $i$ vectors of the basis and $p:\sO^{2n}\rightarrow \sO^{2n}/E_n$ is the projection. Then $(E_1,\dots,E_{2n-1})$ is an element of $\text{\rm Rat}_{\underline{k}}\bigl(F(2n-1)\bigr)$ and we obtained the inverse map.\end{proof}
This proposition applies to $T^\ast GL(n,\cx)$, since we have $k_{n+i}-k_{n+i+1}=1$, $i=0,\dots,n-1$. In this case, hovwever, we have a different characterisation of the relevant rational maps. We adopt the following definition (cf. \cite{Craw}):
\begin{definition} A rational map $f:\oP^1\rightarrow G/P$, where $G$ is reductive and $P$ is parabolic, is said to be {\em full} if the image of $f$ is not contained in  any $G^\prime/P^\prime$, for proper subgroups $G^\prime\leq G$ and  $P^\prime \leq P\cap G^\prime$ ($G^\prime$ - reductive, $P^\prime$ - parabolic $G^\prime$).\end{definition} 
\begin{lemma} The stratum $\text{\rm Rat}_{(1,\dots,n);(1,\dots,1)}\bigl(GL(2n,\cx)/P\bigr)$ of rational maps into the manifold of half-full flags in $\cx^{2n}$ coincides with the set of full rational maps of degree $(1,\dots,n)$\end{lemma}
\begin{proof} Let $f=(E_1,\dots,E_n)$ be a rational map of degree $(1,\dots,n)$ into $GL(2n,\cx)/P$. If $\sO^{2n}/E_n\simeq \bigoplus_{i=1}^n \sO(m_i)$, then each $m_i$ is nonnegative, and, since $\sum_{i=1}^n m_i=n$, we have that either $f\in \text{\rm Rat}_{(1,\dots,n);(1,\dots,1)}\bigl(GL(2n,\cx)/P\bigr)$ or one of the $m_i$, say $m_n$, is equal to zero. In the first case, $f$ is full by Proposition \ref{SSSS}. In the second case, let $F=p^{-1}\bigl(\bigoplus_{i=1}^{n-1} \sO(m_i)\bigr)$, where $p:\sO^{2n}\rightarrow \sO^{2n}/E_n$ is the projection. Then $F\simeq\sO^{2n-1} $ and $(E_1,\dots,E_n)$ is a flag of subspaces of $F$. Therefore $f$ cannot be full.
\end{proof}

We now recall from remarks made after Theorem 7.2 that the map \eqref{pi_star} for $T^\ast GL(n,\cx)$ agrees with the map defined in terms of Nahm equations by restricting a solution, defined on the $2n-1$ intervals, to the first $n$ intervals. In particular, we obtain a symplectic structure on $\sR=\text{\rm Rat}_{(1,\dots,n);(1,\dots,1)}\bigl(GL(2n,\cx)/P\bigr)$ and an action of $A\simeq \cx^{n(n+1)/2}$ on $\sR$, given by the polar part of a rational map, as in Section \ref{A(c)}. We conclude:
\begin{theorem} There exists a canonical $A\times GL(n,\cx)$-equivariant symplectomorphism, preserving the relevant moment maps, between  $T^\ast GL(n,\cx)$ and the space $\text{\rm FRat}_{(1,\dots,n)}\bigl(GL(2n,\cx)/P\bigr)$ of full rational maps of degree $(1,\dots,n)$  into the manifold of half-full flags in $\cx^{2n}$.
\hfill $\Box$
\end{theorem}

\appendix
\section{Lax equation and rational maps}

A theorem of Hurtubise \cite{Hu} identifies spaces of rational maps into full flag manifolds with certain moduli spaces of solutions to Nahm's equations. In this appendix we state this result in a form needed to prove Theorems \ref{Phi} and \ref{Psi}. 

\subsection{Lax equation} Nahm's equations are certain gauge-invariant $\u(n)$-valued ODE's, whose moduli spaces carry a natural hyperk\"ahler structure. Since we are interested only in the complex-symplectic structure, and not the metric, we shall only discuss the complex part of Nahm's equations. This takes the form of the {\em Lax equation}:
\begin{equation} \frac{d\beta}{dt}=[\beta,\alpha],\label{Lax}\end{equation}
where $\alpha(t),\beta(t)$ are $\gl(n,\cx)$-valued functions on an interval. 
The space of solutions is acted upon by the gauge group ${\sG}$ of $GL(n,\cx)$-valued functions $g(t)$:
\begin{eqnarray} \alpha&\mapsto & \Ad(g)\alpha -\dot{g}g^{-1}\nonumber\\ \beta &\mapsto & \Ad(g)\beta \label{action3}\end{eqnarray}
By imposing boundary conditions, we define certain moduli spaces of $\gl(n)$-valued solutions $F_{n}(m,r;a,b)$, which are building blocks for other, more complicated, moduli spaces. Here $m,r$ are  nonnegative integers  less than or equal to $n$ and $a<b$ are real numbers. In what follows, a {\em standard} irreducible $k$-dimensional representation of $\fS\fL(2,\cx)$ refers to the following pair of $k\times k$ matrices:
\begin{equation}\diag\left(-\frac{k-1}{2},-\frac{k-3}{2},\dots,\frac{k-1}{2}\right),\qquad
\left(\begin{array}{ccccl}0 & 0 &\dots & 0 & 0\\ 1 & 0 &\dots & 0 & 0\\0 & 1 &\dots & 0 & 0\\
\vdots &\vdots & \ddots &\vdots &\vdots\\0 & 0 &\dots & 1& 0\end{array}\right).\label{residues}\end{equation}
\begin{itemize}
\item Solutions in $F_n(m,r;a,b)$ are defined on $(a,b)$.
\item  The $m\times m$ (resp. $r\times r$) upper-left minors of a solution $(\alpha,\beta)$ in  $F_n(m,r;a,b)$ are analytic at $t=a$ (resp. at $t=b$), while the $(n-m)\times(n-m)$ (resp. $(n-r)\times(n-r)$ lower-right minors have simple poles with residues defining the standard $(n-m)$-dimensional (resp. $(n-r)$-dimensional) irreducible representation of $\fS\fL(2,\cx)$. The off-diagonal blocks are of the form ${(t-a)}^{(n-m-1)/2}\times(\text{\it analytic in $t-a$})$ (resp. ${(t-b)}^{(n-r-1)/2}\times(\text{\it analytic in $t-b$})$.
\item The gauge group  for $F_n(m,r;a,b)$ consists of gauge transformations $g$ with $g(a)=g(b)=1$, such that the derivatives of the off-diagonal blocks behave at $a$ and $b$ as above.
\end{itemize}

The moduli spaces $F_n(m,r;a,b)$ is a smooth complex manifold \cite{BielLMS} and it has a natural holomorphic symplectic structure given by: 
\begin{equation}\int_a^b\tr  d\alpha\wedge d\beta.\label{dd}\end{equation}
Moreover,  $F_n(m,r;a,b)$
admits a tri-Hamiltonian action of $GL(m,\cx)\times GL(r,\cx)$, given by gauge transformations $g$, whose $m\times m$ (resp. $r\times r$) upper-left minor takes arbitrary values at $t=a$ (resp. $t=b$).  Both $GL(m,\cx)$ and $GL(r,\cx)$ act freely. The moment map for the action of $GL(m,\cx)$ is $\beta(a)_{(m)}$, where the subscript $(m)$ denotes the $m\times m$ upper-left minor. Similarly, the moment map for $GL(r,\cx)$ is $-\beta(b)_{(r)}$. 

\subsection{The structure of $F_n(m,r;a,b)$\label{cs}} 

It is easy \cite{Hu,BielCMP} to identify the complex symplectic structure of $F_n(m,n;a,b)$ and $F_n(n,r;a,b)$, i.e. moduli spaces, which have poles only at one end of the interval. We discuss briefly $F_n(m,n;0,b)$. Suppose first that $m<n$.
Let $e_1,\ldots,e_n$ denote the standard basis of $\cx^n$. There is a unique solution $w_1$ of the equation
\begin{equation} \frac{dw}{dt}=-\alpha w\label{alphato0},\quad w:[0,b]\rightarrow \cx^n,\end{equation}
with
\begin{equation} \lim_{t\rightarrow 0}\left(t^{-(n-m-1)/2}w_1(t)-e_{m+1}\right)=0\label{w1}.\end{equation}
Setting $w_i(t)=\beta^{i-1}(t)w_1(t)$, we obtain a solution to \eqref{alphato0} with 
$$ \lim_{t\rightarrow 0}\left(t^{i-(n-m-1)/2}w_i(t)-e_{m+i}\right)=0.$$
In addition there are solutions $u_1,\ldots,u_m$ to \eqref{alphato0} whose last $n-m$ components vanish to order $(n-m+1)/2$, and which are linearly independent at $t=0$.
The complex gauge transformation $g(t)$ with $g^{-1}=(u_1,\ldots,u_m,w_1,\ldots,w_{n-m})$ makes $\alpha$ identically zero and sends $\beta(t)$ to the constant matrix (cf. \cite{Hu}):

\begin{equation}B=\left(\begin{array}{ccc|cccc} &  &  & 0 &\ldots & 0 & b_1\\
& X & & \vdots & & \vdots & \vdots\\ &  &  & 0 &\ldots & 0 & b_m \\ \hline 
a_1 & \ldots & a_m & 0 &\ldots & 0 & c_1\\
0 &\ldots & 0 & 1 & \ddots & & c_2\\
\vdots & & \vdots & & \ddots & \ddots & \vdots \\0 &\ldots & 0 & 0&\ldots & 1 & c_{n-m} \end{array}\right).\label{betaconst}\end{equation}

The mapping $(\alpha,\beta)\mapsto (g(b),X,a,b,c)$ gives a biholomorphism between $F_n(m,n;0,b)$ and $GL(n,\cx)\times\gl(m,\cx)\times \cx^{n+m}$ \cite{BielLMS}. The action of $GL(n,\cx)$ is given by the right translations, and the action of $GL(m,\cx)$ is given by $$p\cdot\bigl(g(b),X,a,b,c\bigr)=(pg(b),pXp^{-1},ap^{-1},pb,c),$$ where for the first term we embedded $Gl(m,\cx)$ in $Gl(n,\cx)$ as the $m\times m$ upper-left minor. 
\par
Now consider the remaining case $m=n$ (cf. \cite{Kr}). This time, there are no poles at $t=a$ and we can simply solve the equation $\dot{g}=g\alpha$, $g(a)=1$, for a gauge transformation $g:[0,b]\rightarrow GL(n,\cx)$.  This $g$ makes $\beta$ constant, equal to some $X$. The correspondence $(\alpha,\beta)\mapsto (g(b),X)$ sets up an isomorphism:
\begin{equation} F_n(n,n;a,b)\simeq T^\ast GL(n,\cx).\label{isomorphism1}\end{equation}
From above, we also have: 
\begin{equation}  F_n(0,n;a,b)\simeq V_n,\label{isomorphism2}\end{equation}
(recall the definition of $V_n$ from Section \ref{Vn}) and both these isomorphisms are equivariant symplectomorphisms. Observe that $F_n(n,n;a,b)$ consists of $\gl(n,\cx)$-valued solutions regular at both ends of the interval, while $F_n(0,n;a,b)$ consists of $\gl(n,\cx)$-valued  solutions regular at $t=b$ and having poles with residues defining the standard $n$-dimensional irreducible representation of $\fS\fL(2,\cx)$ at $t=a$. All $F_n(m,n;a,b)$ can,  actually, be  obtained as complex symplectic quotients of $F_n(n,n;a,b)$, i.e. of $T^\ast GL(n,\cx)$ by a unipotent group as in Section \ref{Vn} (see also \cite{BielLMS}). 
\par
The structure of $F_n(n,r;a,b)$ is completely analogous and one observes now that $F_n(m,r;a,b)$ is the symplectic quotient of $F_n(m,n;a,c)\times F_n(n,r;c,b)$ by the diagonal $GL(n,\cx)$, for any $a<c<b$. This follows from the description of the actions and moment maps, given at the end of previous section: the moment map for the diagonal action of $GL(n,\cx)$ being zero is equivalent to the $\beta$-components of solutions matching at $t=c$; quotienting then by $GL(n,\cx)$ is equivalent to matching the $\alpha$'s. We can, therefore, describe $F_n(n,r;a,b)$ as follows:
\begin{proposition} $F_n(m,r;a,b)$ is biholomorphic to the space of $(g,B)\in GL(n,\cx)\times \gl(n,\cx)$ such that $B$ is of the form \eqref{betaconst} and $g^{-1}Bg$ is of the form \eqref{betaconst} with $r$ replacing $m$.\label{BB}\end{proposition}

\subsection{Hurtubise's theorem}  The main result of Hurtubise's paper \cite{Hu} can be now interpreted as saying that the symplectic quotient of the product of several  ${F}_{n_i}(m_i,r_i;a_i,b_i)$ (and, for some degrees, with some additional flat factors) is canonically biholomorphic to the space $\Rat$ of rational maps. 
\par
Let  $n$ be a natural number, $\underline{k}=(k_1,\dots,k_n)$ a sequence of nonnegative integers, and let $t_1,\dots,t_{n+1}$ be a strictly increasing sequence of real numbers. For every $j=1,\dots,n$, let $m_j=\min\{k_{j-1},k_j\}$ and $r_{j}=\min\{k_{j},k_{j+1}\}$, with the convention that $k_0=k_{n+1}=0$. In addition, put 
$$ s_j=\begin{cases} k_j & \text{if $k_j=k_{j-1}$}\\ 0& \text{if $k_j\neq k_{j-1}$}.\end{cases}$$
Finally set
$$E_j=F_{k_j}(m_{j},r_{j};t_j,t_{j+1})\times \oH^{s_j}$$
for $j=1,\dots,n$. Observe that, for every $j$, $r_{j-1}=m_j$ and  $E_{j-1}\times E_j$ admits a Hamiltonian action of $GL(m_{j},\cx)$, so that the symplectic quotient will again consist of matching solutions at $t=t_j$. 
\par
The main result of \cite{Hu} can be now restated as follows:
\begin{theorem}[Hurtubise] Let $n$ be a natural number, $\underline{k}=(k_1,\dots,k_n)$ a sequence of nonnegative integers, and let $t_1,\dots,t_{n+1}$ be a strictly increasing sequence of real numbers.
The space $\Rat$ is canonically biholomorphic to the symplectic quotient of 
$\prod_{j=1}^{n}E_j$ by $\prod_{j=1}^{n-1}GL(r_j,\cx)$, where each $GL(r_j,\cx)$ acts by changing values of gauge transformations at $t=t_{j+1}$.\label{Hurt}\end{theorem}

Thus, $\Rat$ can be identified with a moduli space of solutions to \eqref{Lax} on the union of intervals $(t_j,t_{j+1})$, $j=1,\dots,n$, with appropriate matching conditions at the points $t_j$. In particular, the polar part $q_1(z),\dots,q_n(z)$, given by \eqref{q}, is identified with
\begin{equation} q_j(z)=\det\bigl(z-\beta(t)\bigr)\enskip \text{for $t\in(t_j,t_{j+1})$}.\label{plar}\end{equation}
From the above theorem, one can extract a description of $\Rat$ in terms of matrices; see Proposition \ref{oooo}.
\begin{remark} While the choice of the $t_j$ is irrelevant to the complex or complex-symplectic structure, it is relevant to the hyperk\"ahler metric. By keeping track of the $t_j$, we are really thinking  of $\Rat$ as the space of rational maps into the adjoint $U(n+1)$-orbit $\sO(t)$ of $t=\sqrt{-1}\diag(t_1,\dots,t_{n+1})$ together with the canonical isomorphism $\sO(t)\simeq GL(n+1,\cx)/B$. \label{t}\end{remark} 
As mentioned in Section \ref{symplectic}, an extension of Hurtubise's theorem, given in \cite{thesis} (but also essentially proved in \cite{BielCMP}), states that this isomorphism preserves the complex symplectic forms:
\begin{theorem}\cite{thesis,BielCMP} The biholomorphism of Theorem \ref{Hurt} respects the complex symplectic forms, given by \eqref{sympl} on $\Rat_k$, and by \eqref{dd} on  the moduli space of solutions to \eqref{Lax}.\label{Biel}\end{theorem}

Let us, finally, restate Theorem \ref{Hurt}, in the simple cases needed for Theorems \ref{Phi} and \ref{Psi}: \begin{corollary} Let $n\in \oN$ and let $t_1,\dots,t_{n+1}$ be a strictly increasing sequence of real numbers.  The space $\text{\rm Rat}_{1,\dots,n}\bigl(F(n)\bigr)$ is canonically isomorphic to the symplectic quotient of 
$$ \left(\prod_{i=1}^{n-1}F_i(i,i;t_i,t_{i+1})\right)\times F_n(n,0;t_{n},t_{n+1})$$
by the action of $\prod_{i=1}^{n-1} GL(i,\cx)$, where every $GL(i,\cx)$ acts, as described above, by changing values of gauge transformations at $t=t_{i+1}$.\label{Phi1}\end{corollary}
\begin{corollary} Let $n\in \oN$ and let $t_1,\dots,t_{2n}$ be a strictly increasing sequence of real numbers.  The space $\text{\rm Rat}_{1,\dots,n,\dots,1}\bigl(F(2n-1)\bigr)$ is canonically isomorphic to the symplectic quotient of 
$$ \left(\prod_{i=1}^{n}F_i(i,i;t_i,t_{i+1})\right)\times \left(\prod_{i=n+1}^{2n}F_{2n-i}(2n-i,2n-i;t_i,t_{i+1})\right)$$
by the action of $\prod_{i=1}^{n-1} GL(i,\cx)\times\prod_{i=n+1}^{2n}GL(2n-i,\cx)$, where every $GL(i,\cx)$ acts, as described above, by changing values of gauge transformations at $t=t_{i+1}$.\label{Psi1}\end{corollary}
Somewhat more explicitly, $\text{\rm Rat}_{1,\dots,n}\bigl(F(n)\bigr)$ is now interpreted as the moduli  space of solutions on $n$ adjoining intervals, with  solutions $\bigl(\alpha^i(t),\beta^i(t)\bigr)$ on the $i$-th interval being $\gl(i,\cx)$-valued and such that, for $i<n$,  $\bigl(\alpha^{i}(t_{i+1}),\beta^i(t_{i+1})\bigr)$ coincides with  the upper-left $i\times i$ minor of $\bigl(\alpha^{i+1}(t_{i+1}),\beta^{i+1}(t_{i+1})\bigr)$, while at $t_{n+1}$, a solution $\bigl(\alpha^n(t),\beta^n(t)\bigr)$ has poles with residues defining the standard $n$-dimensional irreducible representation of $\fS\fL(2,\cx)$.
\par
Similarly, $\text{Rat}_{1,\dots,n,\dots,1}\bigl(F(2n-1)\bigr)$ is interpreted as the moduli  space of solutions on $2n-1$ adjoining intervals, with solutions being $\gl(i,\cx)$-valued on the first $n$ intervals and $\gl(2n-i,\cx)$-valued on the remaining $n-1$ intervals, again with matching conditions as above.

\begin{ack} This work was done, while the first author was a Humboldt Fellow at the University of G\"ottingen. Both authors thank the Alexander von Humboldt Foundation for their support and  University of G\"ottingen for hospitality.\end{ack}

\end{document}